# Some probability inequalities for multivariate gamma and normal distributions

Thomas Royen

University of applied sciences Bingen, Berlinstrasse 109, D-55411 Bingen, Germany,

e-mail: thomas.royen@t-online.de

## Abstract

The Gaussian correlation inequality for multivariate zero-mean normal probabilities of symmetrical n-rectangles can be considered as an inequality for multivariate gamma distributions (in the sense of Krishnamoorthy and Parthasarathy [5]) with one degree of freedom. Its generalization to all integer degrees of freedom and sufficiently large non-integer "degrees of freedom" was recently proved in [10]. Here, this inequality is partly extended to smaller non-integer degrees of freedom and in particular - in a weaker form - to all infinitely divisible multivariate gamma distributions. A further monotonicity property - sometimes called "more PLOD (positively lower orthant dependent)" - for increasing correlations is proved for multivariate gamma distributions with integer or sufficienty large degrees of freedom.

***Keywords and phrases***: probability inequalities, multivariate normal distribution, multivariate gamma distribution, Gaussian correlation inequality, positive lower orthant dependence, monotonicity properties of multivariate gamma distributions

2010 Mathematics Subject Classification: 60E15

## 1. Introduction and Notation

**Some special notations**: For any $(n \times n)$-matrix $A$ the submatrix of $A$ with row and column indices $i \in M$, $\varnothing \neq M \subseteq \{1,...,n\}$, is denoted by $A_M$. In a similar way $A_{M_1, M_2}$ is the submatrix with row indices $i \in M_1$ and column indices $j \in M_2$. For a matrix with some indices, $A_{indices}^{-1}$ means always $(A_{indices})^{-1}$, e.g. $R_M^{-1} = (R_M)^{-1}$ or $R_{11}^{-1} = (R_{11})^{-1}$. For two $(n \times n)$-matrices $A = (a_{ij}) \geq B = (b_{ij})$ means $a_{ij} \geq b_{ij}$ for all $i, j$, and "$A > B$" is equivalent to "$A \geq B$ and $A \neq B$". A null-matrix is denoted by $O$ and a null-vector simply by 0. The transpose of a matrix $A$ is denoted by $A^t$. The letter $T$ stands always for a pos. def. (or sometimes pos. semi-def.) matrix $T = diag(t_1,...,t_n)$. The probability density function (pdf) of a univariate gamma distribution with the shape parameter $\alpha$ is $g_\alpha(x) = (\Gamma(\alpha))^{-1} \exp(-x) x^{\alpha-1}$, $\alpha, x > 0$, and $G_\alpha(x)$ is the corresponding cumulative distribution function (cdf). The non-central gamma cdf with the non-centrality parameter $y$ is given by

$$G_\alpha(x, y) = \exp(-y) \sum_{k=0}^{\infty} G_{\alpha+k}(x) y^k / k!.$$

Recently, in [10] a simple proof was given for the Gaussian correlation inequality (GCI) in its narrow form for multivariate centered normal probabilities of symmetrical n-rectangles. This inequality asserts for all the events $A_i = \{|X_i| \leq x_i\}, \ i = 1,...,n,$ that

$$P\left(\bigcap_{i=1}^{n} A_i\right) \geq P\left(\bigcap_{i=1}^{n_1} A_i\right) P\left(\bigcap_{i=n_1+1}^{n} A_i\right), \ 1 \leq n_1 < n, \tag{1.1}$$

4where $(X_1,...,X_n)$ has a normal $N_n(0,\Sigma)$ - distribution. W. l. o. g., the covariance matrix $\Sigma$ can be chosen as a correlation matrix $R$. According to [11] the inequality (1.1) implies also the more general form of the GCI:

$P(C_1 \cap C_2) \geq P(C_1)P(C_2)$ for all convex and centrally symmetric sets $C_1, C_2 \subseteq \mathbb{R}^n$.

The inequality (1.1) - with $\Sigma = R$ - can be read as an inequality for the events $B_i = \{|Y_i| \leq y_i = \frac{1}{2}x_i^2\}$, where $(Y_1,...,Y_n)$ has an n-variate gamma distribution (in the sense of Krishnamoorthy and Parthasarathy [5]) with the degree of freedom $\nu = 2\alpha = 1$ and the associated correlation matrix $R$, here denoted as the $\Gamma_n(\alpha, R)$-distribution. It can be defined more generally by its Laplace transform (Lt)

$$|I_n + RT|^{-\alpha}, \tag{1.2}$$

with the identity matrix $I_n$, $T = diag(t_1,...,t_n)$, $t_1,...,t_n \geq 0$, and a set $A_n(R)$ of "admissible" positive values $\alpha$, where admissibility means, that the Laplace inversion of (1.2) yields actually a $\Gamma_n(\alpha, R)$- pdf, which is always true for $2\alpha \in \mathbb{N}$. Unfortunately, $A_n(R)$ is not exactly known for the most correlation matrices $R$. Equivalent conditions for an infinitely divisible (inf. div.) Lt $|I_n + RT|^{-1}$ are found in [1] and [3]. In this case, all positive $\alpha$ are admissible.

Extending the univariate non-central gamma densities to functions $g_\alpha(x, y)$ on $x > 0$, $y \in \mathbb{C}$, some mixture representations can be given for the functions with the Lt (1.2): The non-singular correlation matrix $R$ is called "m-factorial" if

$$R = D + AA^t \Leftrightarrow D^{-1/2}RD^{-1/2} = I_n + BB^t, \tag{1.3}$$

with a pos. def. matrix $D = diag(d_1,...,d_n)$ and an $(n \times m)$-matrix $A$ of the minimal possible rank $m$ and with columns $a_\mu$, which may be real or also pure imaginary. With such an m-factorial $R$ the function

$$g_\alpha(x_1,...,x_n; R) := E\left(\prod_{j=1}^{n} d_j^{-1} g_\alpha(d_j^{-1}x_j, \tfrac{1}{2}b_j S_{2\alpha} b_j^t)\right) \tag{1.4}$$

with the rows $b_j$ of $B$ in (1.3) and an expectation referring to a $W(2\alpha, I_m)$-Wishart (or pseudo-Wishart) matrix $S_{2\alpha}$, $2\alpha \in \mathbb{N}$ or $2\alpha > m-1$, has the Lt (1.2), (see [7] or [8]). With a real matrix $A$ in (1.3) the functions in (1.4) are obviously probability density functions. Their Laplace transforms are obtained with the Lt of the non-central gamma densities followed by integration over $S_{2\alpha}$. With $m = 1$ and a real column $A = a$ all $\alpha > 0$ are admissible. Any non-singular $R$ has always an at most $(n-1)$-factorial representation with $D = \lambda I_n$, where $\lambda$ is the lowest eigenvalue of $R$, and a real matrix $A$ of rank $m \leq n-1$. Hence, all non-integer values $2\alpha > n-2$ are also admissible. In the rest of this paper an "$m$- factorial" correlation matrix $R$ means always an $R$ with a representation as in (1.3), but only with a real matrix $A$.

Let $Q = (q_{ij})$ be the matrix $c^{-1}R^{-1} = c^{-1}(r^{ij})$ with

$$c > \max\{r^{jj} \mid j = 1,...,n\} \text{ and } c > (2\lambda)^{-1}, \tag{1.5}$$

where $\lambda$ is the lowest eigenvalue of $R$. Then $\dot{Q} = Q - I_n$ has a spectral norm $\|\dot{Q}\| < 1$ and only negative elements on its diagonal. According to Griffiths [3], the Lt $|I_n + RT|^{-1}$ (with a non-singular $R$) is inf. div. if and

only if

$$(-1)^k q_{i_1,i_2} q_{i_2,i_3} \cdots q_{i_k,i_1} \geq 0 \tag{1.6}$$

for all subsets $\{i_1,...,i_k\} \subseteq \{1,...,n\}$, $k \geq 3$, ($|I_2 + RT|^{-1}$ is always inf. div.).

Then, with $Z = diag(z_1,...,z_n)$, $z_j = (1+c^{-1}t_j)^{-1}$, we obtain the series expansion (see also sec. 2 in [6])

$$|I_n + RT|^{-\alpha} = |QZ|^\alpha |I_n + \dot{Q}Z|^{-\alpha} =$$

$$|QZ|^\alpha \exp\left(\alpha \sum_{k=1}^\infty \frac{(-1)^k}{k} tr(\dot{Q}Z)^k\right) = \sum_{k_1,...,k_n=0}^\infty q(\alpha;k_1,...,k_n) \prod_{j=1}^n z_j^{\alpha+k_j} \tag{1.7}$$

with coefficients $q(\alpha;k_1,...,k_n)$, which are all non-negative under the conditions (1.5), (1.6). By Laplace inversion we get for all $\alpha > 0$ the $\Gamma_n(\alpha, R)$-cdf

$$G_\alpha(x_1,...,x_n;R) = \sum_{k_1,...,k_n=0}^\infty q(\alpha;k_1,...,k_n) \prod_{j=1}^n G_{\alpha+k_j}(c\, x_j), \tag{1.8}$$

which is a convex combination of products of univariate gamma-cdfs. Only with the condition $c > (2\lambda)^{-1}$ the series in (1.8) is also abs. convergent for any non-singular $R$ and for all $\alpha > 0$, $x_1,...,x_n \leq \infty$, but then, it represents the cdf of a $\Gamma_n(\alpha, R)$-distribution only for admissible $\alpha$-values.

Also with the correlation matrix $Q$, having the elements $q_{ij} = (r^{ii}r^{jj})^{-1/2} r^{ij}$, a series representation

$$G_\alpha(x_1,...,x_n;R) = \sum_{k_1,...,k_n=0}^\infty q(\alpha;k_1,...,k_n) \prod_{j=1}^n G_{\alpha+k_j}(r^{jj}x_j) \tag{1.9}$$

can be derived with only non-negative coefficients under the condition (1.6). Under this condition, the spectral norm of $\dot{Q} = Q - I_n$ is again less than 1, since the assumption of the contrary leads to a contradiction. The series in (1.9) is also abs. convergent for any non-singular $R$, $\alpha > 0$ and $x_1,...,x_n < \infty$.

In [10] actually a more general theorem was proved for $\Gamma_n(\alpha, R)$-distributions with $2\alpha \in \mathbb{N}$ or $2\alpha > n-2$, where the special case $2\alpha = 1$ implies the Gaussian correlation inequality. This theorem and its short proof are reproduced in section 2, since the proof shows the benefits, but also the limitations of the applied method, which serves also as a "model" for the similar proofs of the three following theorems in section 3 an 4. The first aim of the underlying paper is an extension of theorem 1 to smaller non-integer values of $\alpha$, in particular to inf. div. $\Gamma_n(\alpha, R)$-distributions. This will be accomplished partly by theorem 2. Theorem 3 deals with m-factorial $(n \times n)$-correlation matrices with $m \leq n-2$. A further monotonicity property - sometimes called "more PLOD"- for increasing correlations is proved in section 4 for multivariate gamma distributions with integer degrees of freedom $v = 2\alpha$ and at least for all real values $2\alpha > n-2$. In particular with $v = 1$, theorem 4 comprises a generalized form of a conjecture of Šidák in section 3 of [12]. For the multivariate centered normal distribution a global version of the "local" theorem of Bølviken and Joag-Dev in [2] is obtained.

If all $\alpha \geq \tfrac{1}{2}$ would be admissible, then theorem 1 and theorem 4 would also hold for all $\alpha \geq \tfrac{1}{2}$ and, in case of infinitely divisible multivariate gamma distributions, also for all $\alpha > 0$. Theorem 2 and theorem 3 would be dispensable as special cases of theorem 1.

## 2. Proof of the Gaussian correlation inequality extended to some multivariate gamma distributions

**Theorem 1.** Let $R = \begin{pmatrix} R_{11} & R_{12} \\ R_{21} & R_{22} \end{pmatrix}$ be a non-singular correlation matrix with $(n_i \times n_i)$-submatrices $R_{ii}$ and $R_{12} \neq O$. Then, for $R_\tau = \begin{pmatrix} R_{11} & \tau R_{12} \\ \tau R_{21} & R_{22} \end{pmatrix}$, $0 \leq \tau \leq 1$, the cdf $G_\alpha(x_1,...,x_n; R_\tau)$ of the $\Gamma_n(\alpha, R_\tau)$ - distribution is increasing in $\tau$ for all positive numbers $x_1,...,x_n$, $2\alpha \in \mathbb{N}$ and at least for all real values $2\alpha > n-2$. This implies the inequality

$$G_\alpha(x_1,...,x_n; R) > G_\alpha(x_1,...,x_{n_1}; R_{11}) G_\alpha(x_{n_1+1},...,x_n; R_{22}) \geq \prod_{j=1}^n G_\alpha(x_j), \tag{2.1}$$

and in particular, for $2\alpha = 1$, the Gaussian correlation inequality.

**Proof.** All the matrices $R_\tau = (1-\tau)R_0 + \tau R_1$, $0 \leq \tau \leq 1$, are non-singular correlation matrices. With the Lt

$$|I_n + R_\tau T|^{-\alpha} |T|^{-1} \text{ of } G_\alpha(x_1,...,x_n; R_\tau), \; |I_n + R_\tau T| = 1 + \sum_{M \neq \emptyset} |R_{\tau,M}| |T_M|, \; M \subseteq \{1,...,n\},$$

$$R_{\tau,M} = \begin{pmatrix} R_{M_1} & \tau R_{M_1,M_2} \\ \tau R_{M_2,M_1} & R_{M_2} \end{pmatrix}, \; \emptyset \neq M_1 \subseteq \{1,...,n_1\}, \; \emptyset \neq M_2 \subseteq \{n_1+1,...,n\}, \; M = M_1 \cup M_2,$$

$\mathrm{r}_{M_1,M_2} = rank(R_{M_1,M_2})$ and the non-negative eigenvalues $\lambda_{M_1,M_2,i}$ (which are squared canonical correlations) of $R_{M_1}^{-1/2} R_{M_1,M_2} R_{M_2}^{-1} R_{M_2,M_1} R_{M_1}^{-1/2}$, $i = 1,...,\mathrm{r}_{M_1,M_2}$, we obtain the non-negative coefficients

$$c_M(\tau) := -\alpha \frac{\partial}{\partial \tau} |R_{\tau,M}| = -\alpha |R_{M_1}||R_{M_2}| \frac{\partial}{\partial \tau} \prod_{i=1}^{\mathrm{r}_{M_1,M_2}} (1-\tau^2 \lambda_{M_1,M_2,i}) = 2\alpha\tau |R_{\tau,M}| \sum_{i=1}^{\mathrm{r}_{M_1,M_2}} \frac{\lambda_{M_1,M_2,i}}{1-\tau^2 \lambda_{M_1,M_2,i}},$$

$0 < \tau < 1$.

Therefore,

$$H_\alpha^*(t_1,...,t_n; R_\tau) := \frac{\partial}{\partial \tau} |I_n + R_\tau T|^{-\alpha} |T|^{-1} = \sum_{M = M_1 \cup M_2} c_M(\tau) |I_n + R_\tau T|^{-(\alpha+1)} |T_M| |T|^{-1} \text{ is the Lt of}$$

$$H_\alpha(x_1,...,x_n; R_\tau) = \sum_{M = M_1 \cup M_2} c_M(\tau) \left(\prod_{i \in M} \frac{\partial}{\partial x_i}\right) G_{\alpha+1}(x_1,...,x_n; R_\tau),$$

which is always positive for any positive numbers $x_1,...,x_n$, since the coefficients $c_M(\tau)$ cannot identically vanish because of the positive rank of $R_{12}$.

Finally, the identity $H_\alpha(x_1,...,x_n; R_\tau) = \frac{\partial}{\partial \tau} G_\alpha(x_1,...,x_n; R_\tau)$ follows from the integration of $H_\alpha^*(t_1,...,t_n; R_\vartheta)$ over $\vartheta \in [0,\tau]$, using the uniqueness of the Laplace inversion and Fubini's criterion, since

$$\int_{\mathbb{R}_+^n} (G_\alpha(x_1,...,x_n; R_\tau) - G_\alpha(x_1,...,x_n; R_0)) \prod_{j=1}^n e^{-t_j x_j} dx_j = \int_0^\tau H_\alpha^*(t_1,...,t_n; R_\vartheta) d\vartheta =$$

$$\int_{\mathbb{R}_+^n} \left(\int_0^\tau H_\alpha(x_1,...,x_n; R_\vartheta) d\vartheta\right) \prod_{j=1}^n e^{-t_j x_j} dx_j, \text{ which concludes the proof.} \qquad \square$$

**Remarks**. For identical $x_j = x$ and positive correlation means

$$r_1 = \frac{2}{n_1(n_1-1)} \sum_{1 \leq i < j \leq n_1} r_{ij}, \quad r_2 = \frac{2}{n_2(n_2-1)} \sum_{n_1+1 \leq i < j \leq n} r_{ij},$$

the approximation

$$G_\alpha(x,...,x;R) \approx$$

$$G_\alpha(x,...,x;R_{11})G_\alpha(x,...,x;R_{22}) + \sum_{k=1}^{\infty} \frac{\Gamma(\alpha+k)}{\Gamma(\alpha)k!}((r_1 r_2)^{-1} r^2))^k c_k(x;\alpha,n_1,r_1)c_k(x;\alpha,n_2,r_2), \text{ with}$$

$$r^2 = \frac{1}{n_1 n_2} \sum_{i=1}^{n_1} \sum_{j=n_1+1}^{n} r_{ij}^2 \leq r_1 r_2 \text{ (in the original article only } r = \frac{1}{n_1 n_2} \sum_{i=1}^{n_1} \sum_{j=n_1+1}^{n} r_{ij} \text{ ) and}$$

$$c_k(x;\alpha,n_i,r_i) = \frac{\Gamma(\alpha)k!}{\Gamma(\alpha+k)} \int_0^\infty (G_\alpha((1-r_i)^{-1}x,(1-r_i)^{-1}r_i y))^{n_i} L_k^{(\alpha-1)}(y) g_\alpha(y) dy,$$

where the $L_k^{(\alpha-1)}$ are the generalized Laguerre polynomials, is proposed in [9]. This approximation is recommended in particular for small exceedance probabilities $1 - G_\alpha(x,...,x;R)$. The error of this approximation tends to zero with a decreasing variability of the correlations within $R_{11}$, $R_{22}$ and $R_{12}$. Some numerical examples with very accurate conservative approximations of this type are found in [9] and [4].

## 3. Some supplements to theorem 1

For inf. div. $\Gamma_n(\alpha,R)$-distributions a theorem similar to theorem 1 is proved here, but only for a partition of $R$ with $n_1 = n-1$. The inequality (3.1) in theorem 2 means, that such a $\Gamma_n(\alpha,R)$-random vector has positively lower orthant dependent components for all $\alpha > 0$. A simple special case is the following one: For a one-factorial $(n \times n)$-correlation matrix $R \neq I_n$ with correlations $r_{ij} = a_i a_j$, $i \neq j$, and real numbers $a_i \in (-1,1)$ we obtain from (1.4)

$$G_\alpha(x_1,...,x_n;R) = \int_0^\infty \left(\prod_{j=1}^n G_\alpha((1-a_j^2)^{-1}x_j,(1-a_j^2)^{-1}a_j^2 y)\right) g_\alpha(y) dy$$

for all $\alpha > 0$. Since the univariate non-central gamma cdf $G_\alpha(x,y)$ is decreasing in $y$, the inequality

$$G_\alpha(x_1,...,x_n;R) \geq G_\alpha(x_1,...,x_{n_1};R_{11})G_\alpha(x_{n_1+1},...,x_n;R_{22})$$

follows from Kimball's inequality (see e.g. sec. 2.2 in [13]). The inequality is strict for all positive numbers $\alpha$ and $x_1,...,x_n$, if $R_{12} \neq O$.

Within the proof of theorem 2 we use the criterion of Bapat [1] for the infinite divisibility of the Lt $|I_n + RT|^{-1}$ with a non-singular $R$. A signature matrix $S = diag(s_1,..,s_n)$ has only the values $s_i = \pm 1$ on its diagonal. A real non-singular $(n \times n)$-matrix $A$ is an "M-matrix" if $A$ has only non-positive off-diagonal elements and $A^{-1}$ has only non-negative elements. Then, according to Bapat, $|I_n + RT|^{-1}$ is inf. div. if and only if there exists a signature matrix $S$ for which $(SRS)^{-1}$ is an M-matrix.

**Theorem 2**. Let $R = (r_{ij}) = \begin{pmatrix} R_{11} & r \\ r^t & 1 \end{pmatrix}$, $r \neq 0$, be a non-singular $(n \times n)$ - correlation matrix with an inf. div. Lt $|I_n + RT|^{-1}$. Then, for all $R_\tau = \begin{pmatrix} R_{11} & \tau r \\ \tau r^t & 1 \end{pmatrix}$, $0 \leq \tau \leq 1$, the Lt $|I_n + R_\tau T|^{-1}$ is inf. div. and the $\Gamma_n(\alpha, R_\tau)$-cdf $G_\alpha(x_1, ..., x_n; R_\tau)$ is an increasing function of $\tau$ for all positive numbers $\alpha, x_1, ..., x_n$, which implies

$$G_\alpha(x_1, ..., x_n; R) > G_\alpha(x_1, ..., x_{n-1}; R_{11}) G_\alpha(x_n) \geq \prod_{j=1}^n G_\alpha(x_j). \tag{3.1}$$

If $\tau = (\tau_1, ..., \tau_n) \in [0,1]^n$ and $R_\tau = diag(1 - \tau_1^2, ..., 1 - \tau_n^2) + (\tau_i r_{ij} \tau_j)$, then $G_\alpha(x_1, ..., x_n; R_\tau)$ is a non-decreasing function for each $\tau_i$.

**Proof**. According to Bapat [1] there exists a signature matrix $S$ with $SR^{-1}S$ being an M-matrix. With

$$Q = (q_{ij}) = \begin{pmatrix} Q_{11} & q \\ q^t & 1 \end{pmatrix} = SRS \text{ and } q^{nn} = |Q|^{-1}|Q_{11}| = (1 - q^t Q_{11}^{-1} q)^{-1} \text{ we obtain}$$

$$Q^{-1} = \begin{pmatrix} Q_{11}^{-1} + q^{nn} Q_{11}^{-1} q q^t Q_{11}^{-1} & -q^{nn} Q_{11}^{-1} q \\ -q^{nn} q^t Q_{11}^{-1} & q^{nn} \end{pmatrix} \tag{3.2}$$

with only non-positive off-diagonal elements and only non-negative elements in $Q$, and in the same way with

$$Q_\tau = SR_\tau S = \begin{pmatrix} Q_{11} & \tau q \\ \tau q^t & 1 \end{pmatrix} \text{ and } q_\tau^{nn} = |Q_\tau|^{-1}|Q_{11}| = (1 - \tau^2 q^t Q_{11}^{-1} q)^{-1} \leq q^{nn}, \text{ the inverse matrix}$$

$$Q_\tau^{-1} = \begin{pmatrix} Q_{11}^{-1} + \tau^2 q_\tau^{nn} Q_{11}^{-1} q q^t Q_{11}^{-1} & -\tau q_\tau^{nn} Q_{11}^{-1} q \\ -\tau q_\tau^{nn} q^t Q_{11}^{-1} & q_\tau^{nn} \end{pmatrix}.$$

All the off-diagonal elements in the left upper block in (3.2) are non-positive and it is $Q_{11}^{-1} q \geq 0$. Therefore, $Q_\tau^{-1}$ is an M-matrix too and $|I_n + R_\tau T|^{-1}$ is inf. div. for all $\tau \in [0,1]$. Now,

$$\frac{\partial}{\partial \tau} |I_n + R_\tau T|^{-\alpha} |T|^{-1} = -\alpha |I_n + R_\tau T|^{-(\alpha+1)} \sum_M \frac{\partial}{\partial \tau} |R_{\tau,M}| |T_M| |T|^{-1},$$

where the summation extends over all index sets $M = M_1 \cup \{n\}$, $\emptyset \neq M_1 \subseteq \{1, ..., n-1\}$. With

$$c_M(\tau) = -\alpha \frac{\partial}{\partial \tau} |R_{\tau,M}| = 2\alpha \tau |R_{M_1}| r_{M_1}^t R_{M_1}^{-1} r_{M_1} \geq 0,$$

where $r_{M_1}$ has the components $r_{in}$, $i \in M_1$, this leads by Laplace inversion to

$$0 < \sum_M c_M(\tau) \left( \prod_{j \in M} \frac{\partial}{\partial x_j} \right) G_{\alpha+1}(x_1, ..., x_n; R_\tau) = \frac{\partial}{\partial \tau} G_\alpha(x_1, ..., x_n; R_\tau), \ 0 < \tau < 1, \tag{3.3}$$

where the identity is justified as in the proof of theorem 1. Because of $r \neq 0$, the coefficients $c_M(\tau)$ cannot identically vanish. The derivatives in (3.3) are positive for all positive numbers $x_1, ..., x_n$ because of the representation (1.8) of an inf. div. $\Gamma_n(\alpha, R)$- cdf. The remaining parts of the theorem are obvious consequences of (3.3). □

**Remarks**. This proof cannot be extended directly in the same way to a partitioned correlation matrix as in theorem 1 with $n_1 \leq n-2$, since in such cases $|I_n + R_\tau T|^{-1}$ is frequently not inf. div. for all $\tau \in (0,1)$. An example is given by the correlation matrices

$$R_\tau = \begin{pmatrix} 1 & 0.55 & 0.3\tau & 0.36\tau \\ 0.55 & 1 & 0.48\tau & 0.5\tau \\ 0.3\tau & 0.48\tau & 1 & 0.52 \\ 0.36\tau & 0.5\tau & 0.52 & 1 \end{pmatrix}, \quad 0 \leq \tau \leq 1.$$

The matrix $R_1^{-1}$ is an M-matrix, but the element $r_{0.5}^{13}$ in $R_{0.5}^{-1}$ is positive and the remaining elements $r_{0.5}^{ij}, i \neq j$, are negative. Hence, $R_{0.5}^{-1}$ is no M-matrix and $|I_4 + R_{0.5}T|^{-1}$ is not inf. div., which is also recognized from the criterion of Griffiths [3]. However, for $n = 4$, we can also chose $n_1 = n_2 = 2$ since $2(\alpha+1) > n-2 = 2$, but for larger values of $n$ it is not sure if in such cases the functions $G_{\alpha+1}(x_1,...,x_n;R_\tau)$ represent the cdf of a probability distribution for all $\tau \in (0,1)$ and small $\alpha$. Consequently, the positivity of the partial derivatives of $G_{\alpha+1}(x_1,...,x_n;R_\tau)$ in (3.3) cannot always be guaranteed.

The following theorem extends theorem 1 to m-factorial $(n \times n)$ - correlation matrices with $m < n-1$. For the special case with $m = 1$ see the remark on Kimball's inequality at the beginning of this section.

**Theorem 3**. Let $R = \begin{pmatrix} R_{11} & R_{12} \\ R_{21} & R_{22} \end{pmatrix} = D + AA^t$ be a non-singular $(n \times n)$-correlation matrix, $n \geq 4$, with $(n_i \times n_i)$-submatrices $R_{ii}$, $\text{rank}(R_{12}) > 0$, a pos. def. diagonal matrix $D$ and a real $(n \times m)$-matrix $A$ with the minimal possible rank $m$, $1 \leq m \leq n-2$. Furthermore, let $R_{22} = D_B + BB^t$ be a representation of $R_{22}$ with a pos. def. diagonal matrix $D_B$ and a real $(n_2 \times k)$-matrix $B$ with the minimal possible rank $k$, $0 \leq k \leq \min(m, n_2 - 1)$, ($k = 0$ means $R_{22} = I_{n_2}$). Then, for $R_\tau = \begin{pmatrix} R_{11} & \tau R_{12} \\ \tau R_{21} & R_{22} \end{pmatrix}$, $0 \leq \tau \leq 1$, the functions $G_\alpha(x_1,...,x_n;R_\tau)$ are increasing on $\tau \in (0,1)$ for all fixed positive numbers $x_1,...,x_n$ and at least for all non-integer values $2\alpha > \max(0, \min(m+k-3, n-4))$ or $2\alpha \in \mathbb{N}$. This implies for all positive numbers $x_1,...,x_n$ the inequality

$$G_\alpha(x_1,...,x_n;R) > G_\alpha(x_1,...,x_{n_1};R_{11})G_\alpha(x_{n_1+1},...,x_n;R_{22}) \tag{3.4}$$

at least for $2\alpha > \max(m-1, \min(m+k-3, n-4))$ or $2\alpha \in \mathbb{N}$.

**Remarks**. The function $G_\alpha(x_1,...,x_n;R_\tau)$ with the Lt $|I_n + R_\tau T|^{-\alpha}|T|^{-1}$ is not necessarily the cdf of a probability distribution for all $\tau \in (0,1)$. However, within the proof $G_{\alpha+1}(x_1,...,x_n;R_\tau)$ is shown to be a cdf. Formula (1.4) shows $G_\alpha(x_1,...,x_n;R)$ to be the cdf of a $\Gamma_n(\alpha, R)$- distribution if $2\alpha > m-1$ or $2\alpha \in \mathbb{N}$. E.g., for $n = 6$ a certain fraction of the $(6 \times 6)$-correlation matrices $R$ is 3-factorial, since the $(6 \times 3)$-matrices $A$ provide 15 free parameters for the 15 scalar products representing the 15 correlations in $R$. Hence, the inequality (3.4) holds for such $R$ at least for $2\alpha = 1$ and all values $2\alpha \geq 2$. Similar considerations are possible for higher dimensions n.

**Proof**. With $D = D_1 \oplus D_2$, $(n_i \times n_i)$ - diagonal $D_i$, $A^t = (A_1^t | A_2^t)$ and $k \leq n - m - 1$, we obtain $R_\tau = D_\tau +$ $A_\tau A_\tau^t$ with $D_\tau = D_1 \oplus (\tau^2 D_2 + (1-\tau^2)D_B)$ and the $(n \times (m+k))$-matrix $A_\tau = \begin{pmatrix} A_1 & O \\ \tau A_2 & (1-\tau^2)^{1/2}B \end{pmatrix}$.

Hence, in $G_{\alpha+1}(x_1,...,x_n;R_\tau)$ all real values $2(\alpha+1) > m+k-1$ with $\alpha > 0$ (or $2\alpha \in \mathbb{N}$) are admissible to obtain the cdf of a $\Gamma_n(\alpha+1, R_\tau)$-distribution. If $k > n-m-1$ we can use a different representation $R_\tau = \lambda_\tau I_n + A_\tau A_\tau^t$, $0 < \tau < 1$, with the lowest eigenvalue $\lambda_\tau$ of $R_\tau$, being positive since $R_\tau$ is not singular, and with an $n \times (n-1)$-matrix $A_\tau$ of rank $r \leq n-1$. Then, $2(\alpha+1) > n-2 \Leftrightarrow 2\alpha > n-4$ is admissible. Therefore, all the derivatives $\left(\prod_{j \in M} \frac{\partial}{\partial x_j}\right) G_{\alpha+1}(x_1,...,x_n;R_\tau)$, $x_1,...,x_n > 0$, $\emptyset \neq M \subseteq \{1,...,n\}$, are positive and the rest of the proof can be taken from the proof of theorem 1.

## 4. A further monotonicity property of some multivariate gamma-distributions

**Theorem 4**. Let $R = (r_{ij})$ and $R_0 = (r_{0,ij})$ be non-singular $(n \times n)$-correlation matrices, $n \geq 2$, with $R > R_0$, all $r_{0,ij} > 0$ and $r_0^{ij} \leq 0$ for all off-diagonal elements in $R_0^{-1}$. Then, for all positive numbers $x_1,...,x_n$, $2\alpha \in \mathbb{N}$ or $2\alpha > n-2$, the cdf $G_\alpha(x_1,...,x_n; R_0 + \tau(R-R_0))$ is increasing on $\tau \in [0,1]$, which implies the inequality

$$G_\alpha(x_1,...,x_n; R_0) < G_\alpha(x_1,...,x_n; R) \tag{4.1}$$

**Remarks**. In particular for $2\alpha = 1$ the inequality (4.1) provides an inequality for multivariate central normal probabilities of symmetrical n-rectangles. In Šidák [12] a local version of this latter inequality was proved for the special case with identical correlations $r_{0,ij} = r > 0$. In a similar theorem of Bølviken and Joag-Dev [2] for multivariate normal distributions the correlations in $R_0$ can be continuously increased to correlations of an $R$ as long as $R^{-1}$ remains an M-matrix, whereas theorem 4 does not need such a restriction for $R$. Theorem 4 together with the theorem of Bølviken and Joag-Dev - generalized by convolutions to integer values $2\alpha$ - shows that (4.1) remains true under the assumptions $2\alpha \in \mathbb{N}$, $R > R_0$, all $r_{ij} > 0$, all $r_{0,ij} \geq 0$ and $r_0^{ij} < 0$ for all $i \neq j$, since theorem 4 can be applied to the pair $R, R_\varepsilon = (r_{\varepsilon,ij})$, where $r_{\varepsilon,ij} = \varepsilon$ if $r_{0,ij} = 0$ and $r_{\varepsilon,ij} = r_{0,ij}$ otherwise, with a sufficiently small $\varepsilon > 0$. Since the $\Gamma_n(\alpha, R)$-distribution is invariant under the transformations $R \to SRS$ with any signature matrices $S$, theorem 4 holds also for a pair $R_0, R$ of correlation matrices if there exists a signature matrix $S$ for which $SR_0S$ and $SRS$ satisfy the above assumptions.

**Proof**. All the matrices $R_\tau = (1-\tau)R_0 + \tau R = R_0 + \tau Q$, $Q = (q_{ij}) > O$, $0 \leq \tau \leq 1$, are non-singular correlation matrices and $G_\alpha(x_1,...,x_n;R_\tau)$ with the Lt $|I_n + R_\tau T|^{-\alpha}|T|^{-1}$ is the cdf of a $\Gamma_n(\alpha, R_\tau)$-distribution at least for $2\alpha > n-2$ or $2\alpha \in \mathbb{N}$. From

$$|I_n + R_\tau T| = 1 + \sum_{M \neq \emptyset} |R_{\tau,M}||T_M|, \quad M \subseteq \{1,...,n\}, \tag{4.3}$$

and

$$|R_{\tau,M}| = |R_{0,M} + \tau Q_M| = |R_{0,M}||I_m + \tau Q_M R_{0,M}^{-1}| = |R_{0,M}| \prod_{i=1}^m (1+\tau \lambda_{M,i})$$

with the cardinality $m$ of $M$ and the eigenvalues $\lambda_{M,i}$ of $Q_M R_{0,M}^{-1}$, we obtain the coefficients

$$c_M(\tau) := -\alpha \frac{\partial}{\partial \tau}|R_{\tau,M}| = -\alpha|R_{\tau,M}|\sum_{i=1}^m \frac{\lambda_{M,i}}{1+\tau\lambda_{M,i}}, \quad (c_M(\tau) = 0 \text{ if } m=1). \tag{4.4}$$

The sums $\sum_{i=1}^m \frac{\lambda_{M,i}}{1+\tau\lambda_{M,i}}$ are non-increasing in $\tau$ (decreasing, if there is at least one $\lambda_{M,i} \neq 0$), and they are equal to $tr(Q_M R_{0,M}^{-1})$ if $\tau = 0$.

The Lt $|I_n + R_0 T|^{-1}$ is inf. div. according to the criteria in [1] and [3]. Therefore, all the Laplace transforms $|I_m + R_{0,M} T|^{-1}$ of the corresponding marginal distributions are inf. div. too. Then, according to Bapat [1], there exists always a signature matrix $S_M$, which generates an M-matrix $S_M R_{0,M}^{-1} S_M$ entailing $S_M R_{0,M} S_M \geq O$. The condition $r_{0,ij} > 0$ for all $i, j$ implies $S_M = \pm I_m$. Therefore, all the off-diagonal elements in $R_{0,M}^{-1}$ are non-positive and $tr(Q_M R_{0,M}^{-1}) \leq 0$ since $Q_M \geq O$ with only zeros on its diagonal. In particular, it is $tr(Q_M R_{0,M}^{-1}) < 0$ if $q_{ij} > 0$ and $M = \{i, j\}$. Now we obtain from (4.3) and (4.4) the Lt

$$H_\alpha^*(t_1,...,t_n; R_\tau) := \frac{\partial}{\partial \tau} |I_n + R_\tau T|^{-\alpha} |T|^{-1} = \sum_{M \neq \emptyset} c_M(\tau) |I_n + R_\tau T|^{-(\alpha+1)} |T_M| |T|^{-1} \text{ of}$$

$$H_\alpha(x_1,...,x_n; R_\tau) = \sum_{M \neq \emptyset} c_M(\tau) \left( \prod_{i \in M} \frac{\partial}{\partial x_i} \right) G_{\alpha+1}(x_1,...,x_n; R_\tau) > 0$$

with not identically vanishing coefficients $c_M(\tau) \geq 0$. The identity

$$H_\alpha(x_1,...,x_n; R_\tau) = \frac{\partial}{\partial \tau} G_\alpha(x_1,...,x_n; R_\tau)$$

follows as in the proof of theorem 1. □

It would be of interest for some multiple statistical tests to apply the following (hypothetical) inequality

$$G_\alpha(x,...,x; R_0) \leq G_\alpha(x,...,x; R) \tag{4.5}$$

(in particular for small exceedance probabilities $p = 1 - G_\alpha(x,...,x; R)$) with $2\alpha \in \mathbb{N}$, identical values $x_j = x$ and an $(n \times n)$-correlation matrix $R_0, n \geq 3$, with identical correlations $r$, where $r$ is the mean of the correlations $r_{ij}$ from $R$. For a non-negative correlation mean $r$ at least a local version of (4.5) (i.e. for all $R$ sufficiently close to $R_0$) is true under the condition

$$\lambda(\alpha, n, r, x) = c_1 + (n-4)c_2 - (n-3)c_3 > 0 \text{ with}$$

$$c_1 = \int_0^\infty ((\alpha - \tfrac{1}{2}) f_1^2 + \tfrac{1}{2}(2ryf_2 - f_1)^2) F^{n-2} g_\alpha(y) dy, \quad c_2 = r \int_0^\infty y f_1^2 (2ryf_2 - f_1) F^{n-3} g_\alpha(y) dy,$$

$$c_3 = 2r^2 \int_0^\infty y^2 f_1^4 F^{n-4} g_\alpha(y) dy, \text{ where}$$

$$F = G_\alpha((1-r)^{-1}x, (1-r)^{-1}ry), \quad f_k = \frac{\partial^k}{\partial x^k} G_{\alpha+k}((1-r)^{-1}x, (1-r)^{-1}ry), \quad k = 1, 2,$$

(see theorem 3 in [9]). This condition can frequently be verified only by a plot of the integrand. In [9] also a Taylor approximation is given for $G_\alpha(x,...,x; R_0 + H)$ by a Taylor polynomial $T_2(\alpha, n, r, x; H)$ of $2^{nd}$ degree with the deviations $h_{ij} = r_{ij} - r$, useful for small values $h_{ij}$ and larger values of $x$. It is

$$T_2(\alpha, n, r, x; H) = \int_0^\infty F^n g_\alpha(y) dy + (c_1 - c_2) H_2 + (c_3 - c_2) H_4,$$

$$H_2 = \sum_{1 \leq i < j \leq n} h_{ij}^2, \quad H_4 = \sum_{i < j, k < m} h_{ij} h_{km}, \quad \{i, j\} \cap \{k, m\} = \emptyset.$$

In particular, for the application to the probability $P\{\max |Z_j| \leq z\}$ with a normal $N(0, R)$-random vector $Z$ we find with

$$F = \tfrac{1}{2}\left(erf\left(\tfrac{z+\sqrt{r}\,y}{\sqrt{2(1-r)}}\right) + erf\left(\tfrac{z-\sqrt{r}\,y}{\sqrt{2(1-r)}}\right)\right), \quad S = \sinh\left(\tfrac{\sqrt{r}\,yz}{1-r}\right), \quad C = \cosh\left(\tfrac{\sqrt{r}\,yz}{1-r}\right)$$

the coefficients

$$c_1 = \pi^{-3/2}\sqrt{2}(1-r)^{-3}\int_0^\infty \exp\left(-\tfrac{2z^2+(1+r)y^2}{2(1-r)}\right)\left(y\sqrt{r}S - zC\right)^2 F^{n-2}\,dy,$$

$$c_2 = 2\pi^{-2}(1-r)^{-5/2}\int_0^\infty \exp\left(-\tfrac{3z^2+(1+2r)y^2}{2(1-r)}\right)\left(y\sqrt{r}S - zC\right)S^2 F^{n-3}\,dy,$$

$$c_3 = 2^{3/2}\pi^{-5/2}(1-r)^{-2}\int_0^\infty \exp\left(-\tfrac{4z^2+(1+3r)y^2}{2(1-r)}\right)S^4 F^{n-4}\,dy.$$

More general Taylor approximations of $2^{nd}$ degree for $G_\alpha(x,...,x;R)$ are also found in [9] and in sec. A.5 in [4].